\documentclass[12pt]{amsart}

\newtheorem{theorem}{Theorem}[section]
\newtheorem{lemma}[theorem]{Lemma}

\theoremstyle{plain}
\newtheorem{definition}[theorem]{Definition}
\newtheorem{example}[theorem]{Example}
\newtheorem{remark}[theorem]{Remark}

\newtheorem{question}[theorem]{Question}

\newtheorem{setting}[theorem]{Setting}

\theoremstyle{definition}

\theoremstyle{remark}

\numberwithin{equation}{section}

\usepackage{fancyhdr}
\usepackage{amscd}
\usepackage{amsmath}
\usepackage{amsthm}
\usepackage{amssymb}

\usepackage[all]{xy}

\begin{document}

\title[Determinant and divisors]{Determinant functors and deformation of Chow groups}

\author{Sen Yang }

\address{School of Mathematics \\  Southeast University \\
Nanjing, China\\
}
\address{Shing-Tung Yau Center of Southeast University \\ 
Southeast University \\
Nanjing, China\\
}
\email{101012424@seu.edu.cn}

\date{}

\maketitle

\begin{abstract}
Using determinant functor, we describe a natural transformation from local Hilbert functor to K-theoretic cycle groups of codimension one, which were variants of Balmer's tensor triangular Chow groups. This enables us to answers a question of Bloch on infinitesimal deformation of Chow groups of divisors. 
\end{abstract}


\section{Introduction}

Let $Y \subset X$ be a closed irreducible subvariety of codimension $1$, where $X$ is a smooth projective variety over a field $k$ of characteristic zero. The subvariety $Y$ gives an algebraic cycle (class) of the cycle group $Z^{1}(X)$ (and Chow group $CH^{1}(X)$). The main purpose of this paper is to compare deformation of the subvariety $Y$ with that of algebraic cycles (classes). 

Let $Art_{k}$ denote the category of local artinian $k$-algebra with residue field $k$ and let $Set_{*}$ denote the category of pointed sets. The local Hilbert functor $\mathrm{\mathbb{H}ilb}$ is defined to be
\begin{align*}
\mathrm{\mathbb{H}ilb}: \  Art_{k} & \longrightarrow  Set_{*} \\
  \ A & \longrightarrow \mathrm{\mathbb{H}ilb}(A),
\end{align*}
where $\mathrm{\mathbb{H}ilb}(A)$ denotes the set of infinitesimal embedded deformations of $Y$ in $X\times_{\mathrm{Spec}(k)}\mathrm{Spec}(A)$, which is denoted by $X_A$ below.

One useful tool in this study is the following identification
\begin{equation}
CH^{1}(X)=H^{1}(X, K_{1}(O_{X})),
\end{equation}
where $K_{1}(O_{X})=O^{\ast}_{X}$ is the K-theory sheaf associated to the presheaf $U \to K_{1}(O_{X}(U))$ with $U \subset X$ open affine. For $A \in Art_{k}$, one defines a functor (see page 406 of Bloch \cite{Bl3})
\begin{align}
\widetilde{CH}^{1}: & \ A \to H^{1}(X, K_{1}(O_{X_A})),
\end{align}
and considers $\widetilde{CH}^{1}(A)$ as infinitesimal deformation of the Chow group $CH^{1}(X)$.

The following is a special case of a question suggested by Bloch in the introduction of \cite{Bl3}.
\begin{question}  [\cite{Bl3}]\label{q:Bloch's ques}
With notation as above, is there a natural transformation from local Hilbert functor $\mathrm{\mathbb{H}ilb}$ to the functor $\widetilde{CH}^{1}$
\[
\mathrm{\mathbb{H}ilb} \to \widetilde{CH}^{1}?
\]

\end{question}

Guided by this question, this paper is organized as follows. In the second section, we briefly recall determinant functor. In section 3, we use determinant of perfect complexes to describe a natural transformation from local Hilbert functor to Chow group of divisor, and answer Question \ref{q:Bloch's ques} in Theorem \ref{t: transf-Hilb-Chow}.

To conclude the introduction, we remark that the identification (1.1) is a special case of Bloch formula (cf. Bloch \cite{Bl2-Annals}, Kerz \cite{Kerz}, Quillen \cite{Quillen} and Soul\'e \cite{Soule})
\begin{equation*}
CH^{p}(X)=H^{p}(X, K^{M}_{p}(O_{X})),
\end{equation*}
where $p$ is an integer satisfying that $1\leq p \leq \mathrm{dim}(X)$. One may define a functor $\widetilde{CH}^{p}$ in a similar way as in (1.2) and consider a generalization of Question \ref{q:Bloch's ques} asked by Bloch \cite{Bl3}, which was partially answered in \cite{Y5}.  The functor $\widetilde{CH}^{p}$ is very useful in understanding deformations of Chow groups, see Bloch \cite{Bl3}, Bloch-Esnault-Kerz \cite{BEK1, BEK2}, Green-Griffiths \cite{GGChow}, Morrow \cite{Morrow}, Patel-Ravindra \cite{PatelRavi1}, Stienstra \cite{Stien1,Stien2} and others.

\section{Determinant functor}

The determinant of a square matrix is a fundamental concept. In \cite{KM}, Knudsen-Mumford started the functorial study of determinant and constructed a determinant functor on the category of perfect complexes of $O_{X}$-modules. To motivate Knudsen-Mumford's construction, we recall that, for $V$ a finite dimensional vector space, an elementary example of determinant functor is given as
\[
det: V \to (\wedge^{dimV}V, dimV).
\]

This functor can be generalized to free modules of finite ranks and to locally free sheaves of finite ranks on schemes. Let $M$ be a free module of rank $r$ over a ring $R$, one defines a determinant functor as 
\[
det: M \to (\wedge^{r}M, r).
\]
Let $F$ be a locally free sheaves of finite rank on a scheme $X$. On each connected component $U_{i}$ of $X$, $F$ is of constant rank $r_{i}$,  one locally defines a determinant functor on $F$ as 
\begin{equation}
det: F|_{U_{i}} \to (\wedge^{r_{i}}(F|_{U_{i}}), r_{i}).
\end{equation}

Let $F_{\bullet}$ be a strict perfect complex of $O_{X}$-modules
{\small
\[
 \begin{CD}
  0 @>>> F_{n} @>>>  \cdots @>>> F_{1} @>>> F_{0} @>>> 0,
 \end{CD}
\] }where each $F_{i}$ is a free $O_{X}$-module of finite rank. Knudsen-Mumford \cite{KM} defined the determinant of $F_{\bullet}$ (denoted $det(F_{\bullet})$) to be 
\begin{equation}
det(F_{\bullet})=\bigotimes^{n}_{i=0}(det(F_{i}))^{(-1)^{i}},
\end{equation} where each $det(F_{i})$ is defined in (2.1) and $(det(F_{i}))^{-1}$ denotes the dual of $det(F_{i})$ with $i$ an odd integer. They proved that this definition enjoyed a list of nice properties, including that, for a short exact sequence of perfect complexes $0 \to F_{\bullet1} \to F_{\bullet2} \to F_{\bullet3} \to 0$, there exists a natural isomorphism $det(F_{\bullet2})\cong det(F_{\bullet1}) \otimes det(F_{\bullet3})$. This gave a determinant functor on the category of perfect complexes of $O_{X}$-modules.

Deligne \cite{DE} axiomatized Knudsen-Mumford's construction and developed a theory of determinant functor on exact category. This was further generalized by Knudsen \cite{KN}, Breuning \cite{BR}, Muro-Tonks-Witte \cite{MTW} and others.

A Picard groupoid $\mathcal{P}$ is a symmetric monoidal category such that all morphisms are invertible and tensoring with any $x \in \mathcal{P}$ gives an equivalence 
\[
x \otimes : \mathcal{P} \longrightarrow \mathcal{P}.
\]

\begin{example}
For $X$ a scheme, we denote by $\textbf{Pic}^{\mathbb{Z}}(X)$ the category of graded line bundles over $X$. An object of $\textbf{Pic}^{\mathbb{Z}}(X)$ is a pair $(L, \alpha)$ with $L$ an line bundle over $X$ and $\alpha: X \to \mathbb{Z}$ a locally constant map. The symmetric monoidal structure is $(L, \alpha)\otimes (L', \alpha')=(L\otimes L', \alpha+\alpha')$.

\end{example}

The following definition of a determinant functor follows Deligne \cite{DE} and Muro-Tonks-Witte \cite{MTW}.
\begin{definition}  \label{d:det}
Let $\mathcal{C}$ be an exact category and $\mathcal{P}$ a Picard groupoid. A determinant functor $det :\mathcal{C} \to \mathcal{P}$ consists of a functor from the subcategory of isomorphisms (denoted $\mathcal{C}^{iso}$),
\[
det: \ \mathcal{C}^{iso} \longrightarrow \mathcal{P},
\]
together with the data: for any short exact sequence in $\mathcal{C}$
\[
 \delta:  0 \to A \to B \to C \to 0,
\] a morphism
\[
det(\delta): det(C)\otimes det(A) \to det(B),
\]
such that the following must be satisfied:

(1) Naturality: For an isomorphism of two exact sequences $\delta$ and $\delta'$
{\small
\[
\begin{CD}
\delta: \ 0 @>>> A @>>> B @>>> C @>>> 0 \\
@. @VVV @VVV @VVV @. \\
\delta': \ 0 @>>> A' @>>> B' @>>> C' @>>> 0,
\end{CD}
\]
}it follows that $det(\delta)=det(\delta')$. 

(2) Associativity: For an admissible filtration  $A \rightarrowtail B \rightarrowtail C$ which gives four short exact sequences in $\mathcal{C}$ ($A_{1}$, $B_{1}$ and $C_{1}$ denote the cokernel of $A \rightarrowtail B$, $B \rightarrowtail C$ and $A \rightarrowtail C$ respectively)
\begin{align*}
\delta_{1}: 0 \to A \to B \to A_{1} \to 0, \ & \delta_{2}: 0 \to B \to C \to B_{1} \to 0, \\
\delta^{1}_{2}: 0 \to A \to C \to C_{1} \to 0, \ &  \tilde{\delta}: 0 \to A_{1} \to C_{1} \to B_{1} \to 0,
\end{align*}
the following diagram in $\mathcal{P}$ commutes:

\[
 \xymatrix{
&det(C) & \\
det(B_{1})\otimes det(B) \ar[ur]^{det(\delta_{2})} & & det(C_{1})\otimes det(A) \ar[ul]_{det(\delta^{1}_{2})} \\
det(B_{1})\otimes (det(A_{1}) \otimes det(A)) \ar[u]^{1\otimes det(\delta_{1})} \ar[rr]  & & (det(B_{1})\otimes det(A_{1})) \otimes det(A) \ar[u]_{det(\tilde{\delta}) \otimes 1}. \\
 }
 \]

(3) Commutativity: given two objects $A$ and $B$ in $\mathcal{C}$, there are two short exact sequences

\[
\delta_{1}: 0 \to A \to A\oplus B \to B \to 0, \ \delta_{2}: 0 \to B \to A\oplus B \to A \to 0,
\]
and the following triangle commutes:

\[
 \xymatrix{
  & & det(A \oplus B) &  \\
 & det(B)\otimes det(A)      \ar[ur]^{det(\delta_{1})} \ar@{->}[rr] &  & det(A)\otimes det(B) \ar[ul]_{det(\delta_{2})}. 
 }
 \]
\end{definition}

Determinant functors find many applications, for example, see \cite{BBE,BGS,DE,LO,MTW,S2,Wi1}. The following example is used in next section. 
\begin{example} \label{e:deter}
In notation of Setting \ref{s:setting} below, let $\mathcal{C}$ be the exact category of pseudo-coherent $O_{X_A,y}$-modules supported on $y$ and of Tor-dimension $\leq l$ on $O_{X_A,y}$, see Lemma \ref{l:tt} recalled below. Knudsen-Mumford's construction gave a determinant functor on $\mathcal{C}$ in the sense of Definition \ref{d:det}. Concretely, for an element $M \in \mathcal{C}$, $M$ has a resolution
\begin{equation*}
 \begin{CD}
  L_{\bullet} @>>> M @>>> 0,
 \end{CD}
\end{equation*}where $L_{\bullet}$ is a strict perfect complex of $O_{X_A,y}$-modules. The determinant of $M$ is defined to be the determinant of the complex $L_{\bullet}$
\[
det(M)=det(L_{\bullet}),
\]where $det(L_{\bullet})$ is defined in (2.2).

Let $\mathcal{P}$ be Picard category of graded invertible modules on $O_{X_A,y}$. The above construction gives a determinant functor
\begin{align*}
\mathcal{C} & \longrightarrow \mathcal{P} \\
\ \ \  M    &  \longrightarrow \ det(M).
\end{align*}

\end{example}

\begin{remark} \label{r:BGS}
For a vector bundle $F$, the determinant $det(F)$ in (2.1) is defined to be the graded line bundle $(\wedge^{r}F, r)$, where $r$ is an integer. This is in order to handle the sign problem which is not an issue in the study below, so we omit it and consider $\wedge^{r}F$ as the determinant of $F$ below. Cf. Remark 1.2 of Bismut-Gillet-Soul\'e \cite{BGS}.

\end{remark}

\section{K-theoretic cycles}
We introduce K-theoretic cycles in section 3.1 and answer Question \ref{q:Bloch's ques} in section 3.2. The following setting is used in this section. 
\begin{setting} \label{s:setting}
Let $X$ be a nonsingular projective variety over a field $k$ of characteristic zero, with generic point $\eta$. For any $A \in Art_{k}$, we denote by $X_{A}$ the fibre product $X \otimes_{\mathrm{Spec}(k)} \mathrm{Spec}(A)$.

Let $Y \subset X$ be a closed irreducible subvariety of codimension $1$, with generic point $y$. There exists a finite open affine covering $\{U_{i} \}_{i \in I}$ such that $Y \cap U_{i}$ is defined by $f_{i1}$. For $Y' \in \mathrm{\mathbb{H}ilb}(A)$, $Y' \cap U_{i}$ is defined by $f^{A}_{i1} \in O_{X_A}(U_{i})$, which is a lifting of $f_{i1}$.
\end{setting} 

\subsection{K-theoretic cycles}
Keeping Question \ref{q:Bloch's ques} in mind, we want to find a resolution of the K-theory sheaf $K_{1}(O_{X_A})$. This can be done by using Bloch-Ogus-Gabber theorem (cf. \cite{CTHK}). Concretely, it can be proved that (for example, see \cite{DHY} and \cite{Y5}), the Zariski sheafification of the Bloch-Gersten-Quillen sequence
{\footnotesize
\begin{align}
0  \rightarrow K_{1}(O_{X_A,\eta}) \to \bigoplus_{x \in X^{(1)}}K_{0}(O_{X_A,x} \ \mathrm{on} \ x)  \to \cdots,
\end{align}
}is a flasque resolution of the K-theory sheaf $K_{1}(O_{X_A})$. When $A=k$, the sequence (3.1) is of the form
{\footnotesize
\begin{align}
0  \rightarrow K_{1}(O_{X,\eta}) \to \bigoplus_{x \in X^{(1)}}K_{0}(O_{X,x} \ \mathrm{on} \ x)  \to \cdots,
\end{align}
}whose Zariski sheafification is a flasque resolution of the K-theory sheaf $K_{1}(O_{X})$.

Let $q$ be an integer satisfying that $0 \leq q \leq d$, where $d =\mathrm{dim}(X)$. For $x \in X^{(q)}$, let $\overline{K}_{1-q}(O_{X_{A},x} \ \mathrm{on} \ x)$ be the relative K-group, i.e., the kernel of the map (induced by $A \to k$)
\[
K_{1-q}(O_{X_{A},x} \ \mathrm{on} \ x) \to K_{1-q}(O_{X,x} \ \mathrm{on} \ x).
\]
It follows from (3.1) and (3.2) that the Zariski sheafification of 
{\footnotesize
\begin{align}
0  \rightarrow \overline{K}_{1}(O_{X_A,\eta}) \to \bigoplus_{x \in X^{(1)}}\overline{K}_{0}(O_{X_A,x} \ \mathrm{on} \ x)  \to \cdots,
\end{align}
} is a flasque resolution of the sheaf $O_{X}\otimes m_{A}$, which is the kernel of the morphism $K_{1}(O_{X_A}) \to K_{1}(O_{X})$.

\begin{lemma} \label{l: compute-K-relative}
With notations as above, for $x \in X^{(q)}$, there are isomorphisms (between relative K-group and local cohomology)
\[
\overline{K}_{1-q}(O_{X_{A},x} \ \mathrm{on} \ x)=H^{q}_{x}(O_{X}\otimes m_{A}).
\]

\end{lemma}

\begin{proof}
Since K-theory satisfies Zariski descent, we can identify the relative group $\overline{K}_{1-q}(O_{X_{A},x} \ \mathrm{on} \ x)$ with hypercohomogy
\[
\overline{K}_{1-q}(O_{X_{A},x} \ \mathrm{on} \ x)= \mathbb{H}^{-(1-q)}_{x}(O_{X,x}, \overline{K}(O_{X_{A},x})),
\]
where $\overline{K}(O_{X_{A},x})$ is the fiber of map of K-theory spectra $K(O_{X_{A},x}) \to K(O_{X,x})$.

There exists a spectral sequence (cf. Theorem 10.3 of \cite{TT})
{\footnotesize
\begin{equation}
E_{2}^{i, j} = H^{i}_{x}(O_{X,x}, \overline{K}_{j}(O_{X_{A},x}))
 \Longrightarrow \mathbb{H}^{-(1-q)}_{x}(O_{X,x}, \overline{K}(O_{X_{A},x})),
\end{equation}
}where $i-j=-(1-q)$ and $\overline{K}_{j}(O_{X_{A},x})$ is the kernel of $K_{j}(O_{X_{A},x}) \to K_{j}(O_{X,x})$.

Since the Krull dimension of $O_{X,x}$ is $q$, if $i > q$, then the local cohomology $H^{i}_{x}(O_{X,x}, \overline{K}_{j}(O_{X_{A},x}))=0$ for each $j$. This shows that the index $i$ in non-zero terms of the spectral sequence (3.4) satisfies that $0 \leqslant i \leqslant q$. It follows that $j=i+1-q \leqslant 1$. If $j < 1$, then $\overline{K}_{j}(O_{X_{A},x})=0$. Hence, the index $j=i+1-q$ in non-zero terms of the spectral sequence (3.4) can only be $1$, which implies that $i=q$.

In conclusion, the only non-zero term in the spectral sequence (3.4) is $H^{q}_{x}(O_{X,x}, \overline{K}_{1}(O_{X_{A},x}))$. So the spectral sequence (3.4) degenerates and
{\small
\begin{align*}
\mathbb{H}^{-(1-q)}_{x}(O_{X,x}, \overline{K}(O_{X_{A},x})) = H^{q}_{x}(O_{X,x}, \overline{K}_{1}(O_{X_{A},x})) 
= H^{q}_{x}(O_{X}\otimes m_{A}).
\end{align*}
}

\end{proof}  

Using (3.1), (3.2), (3.3) and Lemma \ref{l: compute-K-relative}, one sees that
\begin{theorem} \label{theorem: Main-Diag}
There exists the following commutative diagram in which the Zariski sheafification of each column is a flasque resolution of $O_{X} \otimes m_{A}$,  $K_{1}(O_{X_{A}})$ and $K_{1}(O_{X})$ respectively. 

The maps from the middle column to the right one, denoted $\mathrm{Pr}$, are induced by augmentation $A \to k$. Since the right column is always a direct summand of the middle column, the diagram is split and there exists maps from the middle column to the left one, denoted $\mathrm{P}$

{\tiny
\[
  \begin{CD}
     0 @. 0 @. 0\\
      @VVV @VVV @VVV\\
     O_{X,\eta} \otimes m_{A} @<\mathrm{P}<< K_{1}(O_{X_{A},\eta}) @>\mathrm{Pr}>>  K_{1}(O_{X,\eta}) \\
      @V \partial_{1,X_{A}}^{0,-1}VV @Vd_{1,X_{A}}^{0,-1}VV @Vd_{1,X}^{0,-1}VV\\
        \bigoplus\limits_{x \in X ^{(1)}} H^{1}_{x}(O_{X})\otimes m_{A} @<<< \bigoplus\limits_{x \in X^{(1)}}K_{0}(O_{X_{A},x} \ \mathrm{on} \ x) @>>>  \bigoplus\limits_{x \in X ^{(1)}}K_{0}(O_{X,x} \ \mathrm{on} \ x) \\
     @V \partial_{1,X_{A}}^{1,-1}VV  @Vd_{1,X_{A}}^{1,-1}VV @Vd_{1,X}^{1,-1}VV\\
    \bigoplus\limits_{x \in X ^{(2)}} H^{2}_{x}(O_{X})\otimes m_{A} @<<< \bigoplus\limits_{x \in X^{(2)}}K_{-1}(O_{X_{A},x} \ \mathrm{on} \ x)
      @>>> \bigoplus\limits_{x \in X ^{(2)}}K_{-1}(O_{X,x} \ \mathrm{on} \ x) \\
      @VVV  @VVV @VVV\\
     \vdots @<<< \vdots @>>> \vdots \\ 
     @VVV @VVV @VVV\\
     \bigoplus\limits_{x \in X ^{(d)}} H^{d}_{x}(O_{X})\otimes m_{A} @<<<  \bigoplus\limits_{x \in X^{(d)}}K_{1-d}(O_{X_{A},x} \ \mathrm{on} \ x) @>>> \bigoplus\limits_{x \in X ^{(d)}}K_{1-d}(O_{X,x} \ \mathrm{on} \ x) \\
     @VVV @VVV @VVV\\
     0 @. 0 @. 0.
  \end{CD}
\]
}
\end{theorem}
After tensoring with $\mathbb{Q}$, we proved a generalization of this theorem by using Goodwillie isomorphism (from K-theory to negative cyclic homology), cf. Theorem 2.21 of \cite{Y5}.

\begin{definition}[Definition 3.4 of \cite{Y2}] \label{definition: MilnorKChow}
In notation of Theorem \ref{theorem: Main-Diag}, the first K-theoretic cycle groups of $X$ and $X_A$, denoted $Z_{1}(D^{\mathrm{perf}}(X))$ and $Z_{1}(D^{\mathrm{perf}}(X_A))$ respectively, are defined to be 
{\footnotesize
\[
  Z_{1}(D^{\mathrm{perf}}(X)):= \mathrm{Ker}(d_{1,X}^{1,-1}), \ Z_{1}(D^{\mathrm{perf}}(X_A)):= \mathrm{Ker}(d_{1,X_A}^{1,-1}).
\]
}
The first K-theoretic Chow groups of $X$ and $X_A$, denoted by $CH_{1}(D^{\mathrm{perf}}(X))$ and $CH_{1}(D^{\mathrm{perf}}(X_A))$ respectively, are defined to be 
{\footnotesize
\[
   CH_{1}(D^{\mathrm{perf}}(X)):= \dfrac{\mathrm{Ker}(d_{1,X}^{1,-1})}{\mathrm{Im}(d_{1,X}^{0,-1})}, \ CH_{1}(D^{\mathrm{perf}}(X_A)):= \dfrac{\mathrm{Ker}(d_{1,X_A}^{1,-1})}{\mathrm{Im}(d_{1,X_A}^{0,-1})}.
\]
}
\end{definition}

The above definitions of K-theoretic cycles and Chow groups are based on Balmer's tensor triangular Chow groups \cite{Ba}, see \cite{Klein,Y2,Y4,Y5} for applications. It was proved that (see Theorem 3.16 of \cite{Y2})
\[
 Z_{1}(D^{\mathrm{perf}}(X))= Z^{1}(X), \ CH_{1}(D^{\mathrm{perf}}(X)) = CH^{1}(X).
\]
Moreover, it is immediately from Theorem \ref{theorem: Main-Diag} that
\begin{equation}
CH_{1}(D^{\mathrm{perf}}(X_A))=H^{1}(X, K_{1}(O_{X_{A}})).
\end{equation}

\subsection{Determinant and transformation}

The following lemma was from 5.7 of \cite{TT}.
\begin{lemma} [\cite{TT}] \label{l:tt}
Let $X$ be a scheme with an ample family
of line bundles. Let $i : Y \to X$ be a regular closed immersion ([SGA 6]
VII Section 1) defined by ideal $J$. Suppose $Y$ has codimension $k$ in X.
Then $K(X \ \mathrm{on} \ Y)$ is homotopy equivalent to the Quillen K-theory $K^{Q}(X \ \mathrm{on} \ Y)$
of the exact category of pseudo-coherent $O_{X}$-modules supported
on the subspace $Y$ and of Tor-dimension $\leq k$ on $X$.
\end{lemma}

The notations of Setting \ref{s:setting} is used here. From now on, we identify $K_{0}(O_{X_{A},y} \ \mathrm{on} \ y)$ with $K^{Q}_{0}(O_{X_{A},y} \ \mathrm{on} \ y)$, which is Grothendieck group of the exact category of pseudo-coherent $O_{X_A,y}$-modules supported on the subspace $y$ and of Tor-dimension $\leq 1$. It is obvious that the module $O_{X_{A},y}/(f^{A}_{i1})$ defines an element (still denoted $O_{X_{A},y}/(f^{A}_{i1})$) of $K^{Q}_{0}(O_{X_A,y}\ \mathrm{on} \ y)$.

\begin{definition} \label{Definition:Hilb-Map-K}
With notation as above, one defines a set-theoretic map 
\begin{align*}
\alpha_{A}: \  \mathrm{\mathbb{H}ilb}(A) &  \longrightarrow K^{Q}_{0}(O_{X_{A},y} \ \mathrm{on} \ y) \\
Y^{'} \ &  \longrightarrow   O_{X_{A},y}/(f^{A}_{i1}).    \notag
\end{align*}

\end{definition}

In the diagram in Theorem \ref{theorem: Main-Diag}, $K_{-1}(O_{X,x} \ \mathrm{on} \ x)=K_{-1}(k(x))=0$, this implies that there exists the following diagram (part of the diagram in Theorem \ref{theorem: Main-Diag})
\begin{equation}
  \begin{CD}
     \bigoplus\limits_{x \in X ^{(1)}}H^{1}_{x}(O_{X})\otimes_{k}m_{A}  @<\mathrm{P}<< \bigoplus\limits_{x \in X ^{(1)}} K^{Q}_{0}(O_{X_{A},x} \ \mathrm{on} \ x) \\
      @V \partial_{1,X_{A}}^{1,-1}VV @V d_{1,X_{A}}^{1,-1}VV     \\
      \bigoplus\limits_{x \in X ^{(2)}}H^{2}_{x}(O_{X})\otimes_{k}m_{A}  @<\cong<< \bigoplus\limits_{x \in X ^{(2)}} K^{Q}_{-1}(O_{X_{A},x} \ \mathrm{on} \ x).
  \end{CD}
\end{equation}

We shall prove that $\alpha_{A}(Y')$ is a K-theoretic cycle (in the sense of Definition \ref{definition: MilnorKChow})
\[
\alpha_{A}(Y') \in Z_{1}(D^{\mathrm{Perf}}(X_{A})).
\]
In other words, we shall verify that $d_{1,X_A}^{1,-1} \circ \alpha_{A}(Y')=0$. By using the commutative diagram (3.6) 
, one sees that this is equivalent to proving that $\partial_{1,X_{A}}^{1,-1} \circ \mathrm{P} \circ \alpha_{A}(Y')=0$.

The determinant functor in Example \ref{e:deter} can be used to describe the map $\mathrm{P}$ in the diagram (3.6). Concretely, let $\mathcal{C}$ be the exact category of pseudo-coherent $O_{X_A,y}$-modules supported on $y$ and of Tor-dimension $\leq 1$ on $O_{X_A,y}$. For an element $M \in \mathcal{C}$, $M$ has a resolution
\begin{equation*}
 \begin{CD}
  L_{\bullet} @>>> M @>>> 0.
 \end{CD}
\end{equation*}Here $L_{\bullet}$ has the form
\begin{equation*}
 \begin{CD}
 0 @>>> L_{1} @>d_{1}>> L_{0},
 \end{CD}
\end{equation*}where $L_{1}$ and $L_{0}$ are free $O_{X_A,y}$-modules with basis $e_{1}, \cdots, e_{r_{1}}$ and $g_{1}, \cdots, g_{r_{0}}$ respectively. We denote by $M_{1}$ the matrix of $d_{1}$, which is an $r_{0}\times r_{1}$ matrix of rank $r_{1}$. We may assume that the $r_{1}$ columns and the first $r_{1}$ rows constitute a submatrix of rank $r_{1}$, denoted $\widetilde{M_{1}}$.

One defines the determinant of $M$ to be the determinant of the complex $L_{\bullet}$ as in Example \ref{e:deter}
\begin{align*}
det(M)&=det(L_{\bullet})=(det(L_{1}))^{(-1)} \otimes det(L_{0}) \\
&\cong Hom(det(L_{1}), det(L_{0})) = Hom(\wedge^{r_{1}}L_{1}, \wedge^{r_{0}}L_{0}),  \notag
\end{align*}
where $(det(L_{1}))^{(-1)}$ is the dual of $det(L_{1})$. And this determinant is well-defined on $K^{Q}_{0}(O_{X_{A},y} \ \mathrm{on} \ y)$.

The morphism $d_{1}$ induces a map $\wedge^{r_{1}}d_{1}: \wedge^{r_{1}}L_{1} \to \wedge^{r_{1}}L_{0}$. It is noted that $g_{r_{1}+1} \wedge \cdots \wedge g_{r_{0}}$ is a free generator of $\wedge^{r_{0}-r_{1}}L_{0}$, which gives a map
\begin{align*}
d': \wedge^{r_{1}}L_{0} & \longrightarrow \wedge^{r_{0}}L_{0} \\
\omega  \ \ & \longrightarrow \omega \otimes g_{r_{1}+1} \wedge \cdots \wedge g_{r_{0}}. 
\end{align*}The composition $d' \circ \wedge^{r_{1}}d_{1}: \wedge^{r_{1}}L_{1} \to \wedge^{r_{0}}L_{0}$ is an element of $det(M)$. Under the isomorphism $det(M) \cong O_{X_A,y}$, one sees that  $d' \circ \wedge^{r_{1}}d_{1}$ corresponds to the classical determinant of the $r_{1}\times r_{1}$ matrix $\widetilde{M_{1}}$: 
\begin{align*}
 det(M) & \longrightarrow O_{X_A,y} \\
 d' \circ \wedge^{r_{1}}d_{1} & \longrightarrow |\widetilde{M_{1}}|.
\end{align*}

For $A \in Art_{k}$, $A=k \oplus m_{A}$, so $O_{X_{A},y}=O_{X,y}\otimes_{k} A=O_{X,y} \oplus (O_{X,y}\otimes_{k} m_{A})$. Let 
\[
\rho: O_{X_{A},y} \to O_{X,y}\otimes_{k} m_{A}
\]
be the projection, for $|\widetilde{M_{1}}| \in O_{X_{A},y}$, $\rho(|\widetilde{M_{1}}|) \in O_{X,y}\otimes m_{A}$.

Let $L_{\bullet}$ denote the Koszul resolution of $O_{X,y}/(f_{i1})$, which has the form
\[
 \ 0 \to L_{1} \xrightarrow{f_{i1}} L_{0},
\]
where $L_{1}=L_{0}=O_{X,y}$. The following diagram
{\small
\begin{equation*}
\begin{cases}
 \begin{CD}
   L_{\bullet} @>>> O_{X,y}/(f_{i1})  \\
  L_{1} @>\rho(|\widetilde{M_{1}}|)>> L_{0} \otimes O_{X,y}\otimes m_{A} (\cong O_{X,y}\otimes m_{A}), 
 \end{CD}
\end{cases}
\end{equation*}
} gives an element $\beta^{A}$ in {\small $Ext^{1}(O_{X,y}/(f_{i1}), O_{X}\otimes m_{A})$}. Since
\[
H_{y}^{1}(O_{X}\otimes m_{A})=\varinjlim_{n \to \infty}Ext^{1}(O_{X,y}/(f_{i1})^{n}, O_{X}\otimes m_{A}),
\]the image $[\beta^{A}]$ of $\beta^{A}$ under the limit is in $H_{y}^{1}(O_{X}\otimes m_{A})$ and it is $\mathrm{P}(M)$. To summerize, 

\begin{lemma} \label{l:map P}
With notations as above, the map $\mathrm{P}$ in diagram (3.6) can be described as 
\begin{align*}
 \mathrm{P}:  K^{Q}_{0}(O_{X_{A},y} \ \mathrm{on} \ y)  & \longrightarrow  H^{1}_{y}(O_{X}\otimes m_{A}) \\
M  & \longrightarrow [\beta^{A}].
\end{align*}

\end{lemma}

In particular, for $\alpha_{A}(Y')=O_{X_A,y}/(f^{A}_{i1}) \in K^{Q}_{0}(O_{X_{A},y} \ \mathrm{on} \ y)$ which has a resolution
\[
 \ 0 \to O_{X_A,y} \xrightarrow{f^{A}_{i1}} O_{X_A,y},
\] the image $\mathrm{P} \circ \alpha_{A}(Y') \in H_{y}^{1}(O_{X}\otimes m_{A})$ is represented by the following diagram (still denoted $\beta^{A}$)
{\small
\begin{equation}
\begin{cases}
 \begin{CD}
   L_{\bullet} @>>> O_{X,y}/(f_{i1})  \\
  L_{1} @>f^{A}_{i1}-f_{i1}>> L_{0} \otimes O_{X,y}\otimes m_{A} (\cong O_{X,y}\otimes m_{A}).
 \end{CD}
\end{cases}
\end{equation}
} 

\begin{remark}
For $X$ of dimension one and $A=k[\varepsilon]/(\varepsilon^{2})$ the ring of dual numbers, the above description agrees with Green-Griffiths' in chapter two of \cite{GGtangentspace}.
\end{remark}

\begin{lemma} \label{l:kernel-HH}
With notation as above, for $\mathrm{P} \circ \alpha_{A}(Y') \in H_{y}^{1}(O_{X}\otimes m_{A})$, one has
\[
\partial_{1,X_{A}}^{1,-1} \circ \mathrm{P} \circ \alpha(Y')=0,
\]
where $\partial_{1,X_{A}}^{1,-1}$ is the differential of the left column of diagram (3.6). 
\end{lemma}

This can be proved by mimicking the argument of Lemma 4.5 of \cite{Y1}. For readers' convenience, we sketch the proof briefly.

\begin{proof}
In notation of Setting \ref{s:setting}, by shrinking $U_{i}$, we assume that $O_{X}(U_{i})$ is local. We extend $f_{i1}$ to a regular system of parameter $\{f_{i1}, f_{i2}, \cdots, f_{id}\}$ of the regular local ring $O_{X}(U_{i})$. The prime ideals $I_{j}:=(f_{i1}, f_{ij})$, where $j=2, \cdots, d$, define generic points $z_{j} \in X^{(2)}$. In the following, to check $\partial_{1,X_{A}}^{1,-1} \circ \mathrm{P} \circ \alpha(Y')=0$, we consider the prime $I_{2}=(f_{i1}, f_{i2})$ which defines the generic point $z_{2}$, other cases work similarly.

Let $I=(f_{i1})$ be the prime ideal which defines the generic point (of $Y$) $y \in X^{(1)}$, then $O_{X,y}=(O_{X,z_{2}})_{I}$. Then $\beta^{A}$ (cf. (3.7)), which represents $\mathrm{P} \circ \alpha(Y')$, can be rewritten as
{\footnotesize
\begin{equation*}
\begin{cases}
 \begin{CD}
   F_{\bullet}(f_{i1}) @>>> (O_{X,z_{2}})_{I}/(f_{i1}) \\
  F_{1} @> \dfrac{f_{i2}}{f_{i2}}(f^{A}_{i1}-f_{i1})>> F_{0} \otimes (O_{X,z_{2}})_{I} \otimes m_{A}.
 \end{CD}
 \end{cases}
\end{equation*}
}Here $F_{\bullet}(f_{i1})$ is of the form
{\footnotesize
\[
 \begin{CD}
  0 @>>> F_{1} @>>> F_{0},
 \end{CD}
\]
}where each $F_{1}=F_{0}=(O_{X,z_{2}})_{I}$. Since $f_{i2} \notin I=(f_{i1})$, $f_{i2}^{-1}$ exists in $(O_{X,z_{2}})_{I}$, we can write $f^{A}_{i1}-f_{i1}=\dfrac{f_{i2}}{f_{i2}}(f^{A}_{i1}-f_{i1})$.

The image $\partial_{1,X_{A}}^{1,-1}(\beta^{A})$ is represented by the following diagram (denoted $\gamma$)
{\footnotesize
\[
\begin{cases}
 \begin{CD}
   F_{\bullet}(f_{i1}, f_{i2}) @>>> O_{X, z_{2}}/(f_{i1}, f_{i2}) \\
  F_{2}@>f_{i2}(f^{A}_{i1}-f_{i1})>> F_{0} \otimes O_{X,z_{2}} \otimes m_{A},
 \end{CD}
 \end{cases}
\]
}where the complex $F_{\bullet}(f_{i1}, f_{i2})$ is of the form
\[
 \begin{CD}
  0 @>>>   F_{2} @>M_{2}>>  F_{1} @>>>  F_{0},
 \end{CD}
\]with $F_{2}=\bigwedge^{2}(O_{X, z_{2}}^{ \ \oplus 2})$, $F_{1}=\bigwedge^{1}(O_{X, z_{2}}^{ \ \oplus 2})$ and $F_{0}=O_{X, z_{2}}$.

Let $\{e_{1}, e_{2} \}$ be a basis of $O_{X, z_{2}}^{\ \oplus 2}$, the map $M_{2}$ is 
{\footnotesize
\[
 e_{1}\wedge e_{2}  \longrightarrow f_{i2}e_{1}-f_{i1}e_{2}.
\]
}Since $f_{i2}$ appears in $M_{2}$, one has that
{\footnotesize
\[
\gamma=0 \in Ext_{O_{X,z_{2}}}^{2}(O_{X, z_{2}}/(f_{i1}, f_{i2}), O_{X,z_{2}} \otimes m_{A}).
\]
}Hence, $\partial_{1,X_{A}}^{1,-1}(\beta^{A})=0$.

\end{proof}

The commutativity of diagram (3.6) yields that $d_{1,X_{A}}^{1,-1} \circ  \alpha(Y')=0$. That is,
\begin{theorem} \label{t:AlphaCycle}
In Setting \ref{s:setting}, for $Y' \in \mathrm{\mathbb{H}ilb(A)}$, $\alpha_{A}(Y')$ is a K-theoretic cycle
\[
\alpha_{A}(Y') \in Z_{1}(D^{\mathrm{Perf}}(X_{A})).
\]
\end{theorem}

The K-theoretic cycle $\alpha_{A}(Y')$ defines an element of K-theoretic Chow group (defined in Definition \ref{definition: MilnorKChow}), which further gives an element of $H^{1}(X,K_{1}(O_{X_A}))$ by identification (3.5), denoted $[\alpha_{A}(Y')]$. There is a set-theoretic map
\begin{align}
 \mathrm{\mathbb{H}ilb}(A) & \to  \widetilde{CH}^{1}(A) \\  
  \ Y' & \to [\alpha_{A}(Y')],   \notag
\end{align}
where $\widetilde{CH}^{1}(A)=H^{1}(X,K_{1}(O_{X_A}))$, see (1.2) on page 1. 

Let $f: B \to A$ be a morphism in the category $Art_{k}$, there exists a commutative diagram of sets ( which can be straightforwardly checked)
\begin{equation*}
\begin{CD}
\mathrm{\mathbb{H}ilb}(B) @>(3.8)>>   \widetilde{CH}^{1}(B) \\
@VVV @VVV \\
\mathrm{\mathbb{H}ilb}(A) @>(3.8)>> \widetilde{CH}^{1}(A).  \\
\end{CD}
\end{equation*}
This implies the following result, which answers Question \ref{q:Bloch's ques}.
\begin{theorem} \label{t: transf-Hilb-Chow}
There exists a natural transformation between functors on $Art_{k}$
\begin{equation*}
\widetilde{\mathbb{T}}:  \mathrm{\mathbb{H}ilb} \to \widetilde{CH}^{1},
\end{equation*}
which is defined to be, for any $ A \in Art_{k}$, $\widetilde{\mathbb{T}}(A)$ is (3.8).

\end{theorem}

\textbf{Acknowledgements}.
This paper is a follow-up to \cite{DHY, Y5}, the author thanks Benjamin Dribus, Jerome William Hoffman and Marco Schlichting for many discussions. He also thanks Spencer Bloch, Bangming Deng, Phillip Griffiths, Kefeng Liu and Christophe Soul\'e for their suggestions and/or comments on related work \cite{Y1,Y2,Y4,Y5}.

\end{document}